%% file: main.tex
\documentclass[letterpaper, 10 pt]{ieeeconf}  %

\IEEEoverridecommandlockouts
\usepackage{graphics} %
\usepackage{epsfig} %
\usepackage{amsmath} %
\usepackage{amssymb}  %
\usepackage{booktabs}
\usepackage{makecell}
\usepackage{hyperref}
\usepackage{tikz}

\usepackage{amsfonts}
\usepackage{xcolor}
\usepackage{booktabs}

\usepackage[short]{optidef} %
\usepackage{relsize}

\let\labelindent\relax
\usepackage{placeins}
\usepackage{comment}
\DeclareMathSymbol{\shortminus}{\mathbin}{AMSa}{"39}

\input{oj_commands}
\usepackage[shortlabels]{enumitem}
\usepackage{amsthm} %

\newlength\figureheight
\newlength\figurewidth

\newcommand{\xpred}{x^{k+1}_{\mathrm{pred}}}
\newcommand{\lin}{\mathrm{lin}}
\newcommand{\A}{\mathrm{A}}
\newcommand{\B}{\mathrm{B}}
\newcommand{\C}{\mathrm{C}}
\newcommand{\D}{\mathrm{D}}

\makeatletter
\def\thm@space@setup{\thm@preskip=2pt
\thm@postskip=2pt}
\makeatother
\newtheoremstyle{newstyle}
{} %
{} %
{\mdseries} %
{} %
{\bfseries} %
{.} %
{ } %
{} %

\theoremstyle{newstyle}

\newtheorem{theorem}{Theorem}

\newtheorem{remark}[theorem]{Remark}

\newtheorem{assumption}[theorem]{Assumption}
\newtheorem{proposition}[theorem]{Proposition}

\title{\LARGE \bf Advanced-Step Real-Time Iterations with Four Levels -- New Error Bounds and Fast Implementation in \acados}

\author{Jonathan Frey$^{1,2}$, Armin Nurkanovi{\'c}$^1$, Moritz Diehl$^{1,2}$%
\thanks{$^{1}$Department of Microsystems Engineering (IMTEK), University Freiburg, Germany,
jonathan.frey@imtek.uni-freiburg.de
}%
 \thanks{$^{2}$Department of Mathematics, University Freiburg, Germany}%
 \thanks{This research was supported by DFG via Research Unit FOR 2401, project 424107692 and 525018088, by BMWK via 03EI4057A and 03EN3054B, and by the EU via ELO-X 953348.}
}

\begin{document}

\maketitle
\thispagestyle{empty}
\pagestyle{empty}

\begin{abstract}

The Real-Time Iteration (RTI) is an online nonlinear model predictive control algorithm that performs a single Sequential Quadratic Programming (SQP) per sampling time.
The algorithm is split into a preparation and a feedback phase, where the latter one performs as little computations as possible solving a single prepared quadratic program.
To further improve the accuracy of this method, the Advanced-Step RTI (AS-RTI) performs additional Multi-Level Iterations (MLI) in the preparation phase, such as inexact or zero-order SQP iterations on a problem with a predicted state estimate.
This paper extends and streamlines the existing local convergence analysis of AS-RTI, such as analyzing MLI of level A and B for the first time, and significantly simplifying the proofs for levels C and D.
Moreover, this paper provides an efficient open-source implementation in \acados, making it widely accessible to practitioners.
\end{abstract}
\vspace{-0.3cm}

\section{Introduction}
Nonlinear Model Predictive Control (NMPC) requires at every sampling instant an approximate online solution of a discrete-time optimal control problem (OCP) of the form
\begin{mini!}
	{\substack{s_0,\ldots, s_N, \\ u_0,\ldots, u_{N \shortminus 1} }}
	{\sum_{i=0}^{N\shortminus 1} L_i(s_i, u_i)  + E(s_N)}
	{\label{eq:acados_OCP}}
	{}
	\addConstraint{s_0}{= x}
	\addConstraint{s_{i+1}\label{eq:acados_OCP_eq}}{=\phi_i(s_i, u_i),}{\ k=0,\ldots,N\!\shortminus 1}
	\addConstraint{0}{\leq h_i(s_i, u_i), \label{eq:acados_OCP_ineq}}{\ k=0,\ldots,N\!\shortminus 1}
	\addConstraint{0}{\leq h_N(s_N)}.
\end{mini!}
Its optimization variables are the states  $s_i \in \R^{n_x} $ at $\tau^i$, $ k = 0, \dots, N $ and the control inputs $ u_i \in \R^{n_u}$ acting on shooting intervals $[\tau^i, \tau^{i+1}]$, $ k = 0, \dots, N\shortminus 1 $.
The values $ s_i $ and $ s_{i+1} $ are coupled by the discrete-time dynamics $\phi_i $, which represent the evolution of the real system over a shooting interval.
The cost is given by the path cost terms $L_i$ %
and the terminal cost term $E$.
The path constraints are given by $h_i$
and the terminal constraint by $h_N$. %
In NMPC, once the initial state $x\in\R^{n_x}$ is known, the parametric nonlinear program (NLP) in \eqref{eq:acados_OCP} is solved, and the first control input $u_0$ is fed back to the plant.

Due to the availability of real-time algorithms and efficient open-source software implementations, NMPC is increasingly used in industrial applications.
Real-time algorithms minimize the feedback delay by reducing the online computational load.
The computations are typically divided into two phases:
1) a preparation phase, in which all computations that can be done without knowing $x$ are performed; and
2) the feedback phase, in which the new control input $u_0$ is computed once $x$ is known.

The Real-Time Iteration (RTI)~\cite{Diehl2001} performs a single Sequential Quadratic Programming (SQP) iteration per sampling time.
In the preparation phase, all function and derivative evaluations necessary to construct a Quadratic Program (QP) are performed.
In the feedback phase, only a single QP is solved.
The Multi-level Iteration (MLI) was introduced in~\cite{Bock2007} and extended in~\cite{Nurkanovic2021b,Wirsching2018}.
It is an SQP-based method that offers several variants, also called levels, which only partially update the QP data to reduce the computation time.
The levels are sorted by the amount of QP data that is updated, where level A only updates the initial state and level D updates all QP data.
The advanced-step controller (ASC)~\cite{Zavala2009} solves in the preparation phase an advanced problem to convergence, i.e., the OCP \eqref{eq:acados_OCP} with a predicted state ${x}_{\mathrm{pred}}$.
In the feedback phase, depending on how active-set changes are handled, it solves a linear system, linear program, or QP~\cite{Jaeschke2014}.
The AS-RTI method~\cite{Nurkanovic2019a,Nurkanovic2020b} combines the RTI and ASC approaches.
In the preparation phase, it approximately solves an advanced problem with some MLI variant, and in the feedback phase, like RTI, it solves a single QP.
This method adds some flexibility by allowing one to trade off computational complexity for numerical errors.
It has good theoretical properties~\cite{Nurkanovic2019a}, but has only been used in prototypical simulation experiments so far~\cite{Nurkanovic2021b,Nurkanovic2020a}.

The open-source \acados{} software package implements efficient algorithms for embedded optimal control, with a focus on SQP-type algorithms that thoroughly exploit the block structure of optimal control problems, such as \eqref{eq:acados_OCP}.
It is written in C and relies on the high-performance linear algebra package BLASFEO~\cite{Frison2018}.
The high-level interfaces to MATLAB and Python, its flexible problem formulation, and the variety of solver options have made it an attractive option for real-world NMPC applications.

\textit{Contributions.}
This paper extends the previously existing analysis of AS-RTI~\cite{Nurkanovic2019a,Nurkanovic2020b}, such as analyzing Level B iterations for the first time, which are computationally cheap and converge to suboptimal but feasible solutions of the original problem~\cite{Bock2007}.
It streamlines and significantly simplifies the proofs for the other two levels, C, D, compared to the previous AS-RTI papers~\cite{Nurkanovic2019a,Nurkanovic2020b}.
Moreover, it presents an efficient implementation of the AS-RTI method in \acados~\cite{Verschueren2021}.
Implementation details and the various algorithmic options are discussed.
We extensively test different variants of the implemented algorithm on a benchmark example and demonstrate in a Pareto plot how one can trade-off computational complexity for optimality.

\textit{Outline.}
Section~\ref{sec:as_rti} recalls the AS-RTI method, Section~\ref{sec:theory} derives novel error bounds, and Section~\ref{sec:implementation} discusses implementation details.
Section~\ref{sec:experiments} provides extensive numerical experiments and Section~\ref{sec:conclusion} summarizes the paper.
\vspace{-.3cm}
\section{The Advanced-Step Real-Time Iteration}\label{sec:as_rti}
This section presents AS-RTI and algorithmic ingredients.
\vspace{-.6cm}
\subsection{Sequential Quadratic Programming} \label{sec:sqp}
The parametric NLP \eqref{eq:acados_OCP} can be written more compactly~as
\vspace{-0.4cm}
\begin{mini!}
	{\substack{w\in\R^{n_w}}}
	{f(w)}
	{\label{eq:parametric_nlp}}
	{}
	\addConstraint{0}{= g(w) + Mx \label{eq:parametric_nlp_eq}}
	\addConstraint{0}{\leq h(w), \label{eq:parametric_nlp_ineq}}
\end{mini!}
where the parameter $x$ enters the equality constraint linearly and $M$ is an embedding matrix of appropriate size.
Parametric NLPs can always be brought into this form, by introducing auxiliary variables and a linear equality constraint~\cite{Diehl2001}. %
The functions $ f:\R^{n_w} \! \to \! \R$, $g: \R^{n_w}\! \to \! \R^{n_g}$, $h: \R^{n_w}\! \to \! \R^{n_h}$
are assumed to be twice continuously differentiable.

The Langrangian function of \eqref{eq:parametric_nlp} is $\mathcal{L}(z)  = f(w) - \lambda^{\top} (g(w)+Mx) - \mu^{\top} h(w)$ with $z=(w,\lambda,\mu)$ and Lagrange multipliers $\lambda \in \R^{n_g}, \mu \in \R_{\geq0}^{n_h}$ corresponding to~\eqref{eq:parametric_nlp_eq} and~\eqref{eq:parametric_nlp_ineq}, respectively.
We denote the primal-dual solution of~\eqref{eq:parametric_nlp} by~$\bar{z}(x)=(\bar{w}(x),\bar{\lambda}(x),\bar{\mu}(x))$, which is under suitable assumptions locally unique~\cite{Nocedal2006}.

In the following, we consider a sequence of parameters $\{x^k\}_{k\geq0}$, i.e., the state at each sampling time $t^k$, and perform a fixed number of iterations for a fixed parameter $x^k$.
We use the following notation.
The solution $\bar{z}(x^k)$ for a fixed $x^k$ is often abbreviated to $\bar{z}^k$.
Furthermore, by performing $j$ iterations for a fixed $x^k$, we compute an approximation $z^{k,j}\approx \bar{z}^k$.
If only a single iteration is performed, as in RTI, index $j$ is omitted and we write $z^k = z^{k,1}$.

For a fixed parameter $x^k$, a local minimizer of NLP~\eqref{eq:parametric_nlp} can be computed by an SQP-type method~\cite{Nocedal2006}.
Given a solution guess $z^{k,0}= (w^{k,0},\lambda^{k,0},\mu^{k,0})$ sufficiently close to a local minimizer, a sequence of QPs is solved
\vspace{-.2cm}
\begin{mini!}[2]
	{\substack{\Delta w}}
	{(a^{k,j})^\top \Delta w + \tfrac{1}{2} \Delta w^{\top} A^{k,j} \Delta w }
	{\label{eq:qp_sqp}}{}
	\addConstraint{g^{k,j}+ M x^k+ G^{k,j} \Delta w}{=0\label{eq:qp_sqp_eq}}
	\addConstraint{h^{k,j}+ H^{k,j} \Delta w }{\geq  0\label{eq:qp_sqp_ineq}}.
\end{mini!}
The symmetric positive definite matrix $A^{k,j} $ %
is an approximation of the Lagrange Hessian $\nabla_{ww}^2 \mathcal{L}(z^{k,j})$.
The vector $a^{k,j} = \nabla_w f(w^{k,j})$ denotes the objective gradient, $g^{k,j} =g(w^{k,j})$, $h^{k,j} =h(w^{k,j})$ the constraint residuals, and $G^{k,j} = \nabla_w g(w^{k,j})^\top$, $H^{k,j} = \nabla_w h(w^{k,j})^\top$ are the Jacobians of the constraints.
The primal-dual solution of~\eqref{eq:qp_sqp} is denoted $(\Delta w^{k,j},\lambda^{k,j}_{\mathrm{QP}},\mu^{k,j}_{\mathrm{QP}})$, and a full SQP step updates the iterates by setting $w^{k,j+1} = w^{k,j}  + \Delta w^{k,j}, \ \lambda^{k,j+1} = \lambda^k_{\text{QP}}$, and $\mu^{k,j+1} = \mu^{k,j+1}_{\text{QP}}$.
Every fixed parameter $x^k$ results in a different sequence of iterates and corresponding QP data, indexed by $j$, hence the use of both indices.
\subsection{Real-time NMPC algorithms}
The QP that computes the feedback for the new parameter $x^{k+1}$ at $t^{k+1}$ is prepared in the preparation phase during $[t^k,t^{k+1}]$.
In the RTI, the QP is constructed at the previous output $z^k=z^{k,1}$, computed at time $t^k$.
In AS-RTI, the point~$z^k$ is further refined by computing an approximate solution to an advanced problem - an NLP~\eqref{eq:parametric_nlp} with a predicted parameter $\xpred$.
The improved linearization point~$z^k$ is now denoted by~$z^k_\lin$.
The way the approximation is computed provides a lot of flexibility in algorithm choice and allows for a whole family of different algorithms.
Finally, in the feedback phase, QP~\eqref{eq:qp_sqp} evaluated at $z^k_{\mathrm{lin}} \approx \bar{z}(x^{k+1}_{\mathrm{pred}})$ and $x^{k+1}$ is solved, resulting in the new output $z^{k+1}=z^{k+1,1}\approx \bar{z}(x^{k+1})$.
The main steps of the AS-RTI are:
\begin{enumerate}[(S1)]
	\item At time $t=t^k$: Predict the initial state $x^{k+1}_{\mathrm{pred}}$ at $t^{k+1}$ \label{step:advance}
	\item At $t\in \left[t^k,t^{k+1}\right)$: Starting from the last output $z^k$, iterate on \eqref{eq:parametric_nlp} with $x = x^{k+1}_{\mathrm{pred}}$ with some MLI variant, see Sec.~\ref{sec:mli}, to obtain $z^k_\lin$ -- (``\textit{the inner iterations}''). \label{step:preparation}
	\item At $t\in \left[t^k,t^{k+1}\right)$: Construct QP~\eqref{eq:qp_sqp} on the linearization point $z^k_{\mathrm{lin}}$.
	\item At time $t^{k+1}$, solve \eqref{eq:qp_sqp} with $x = x^{k+1}$. \label{step:feedback}
\end{enumerate}
Note that in the RTI scheme, step \ref{step:preparation} simplifies to setting $z^{k}_{\mathrm{lin}} = z^k$.
In the ASC, $z^{k}_{\mathrm{lin}}$ is a local minimizer of \eqref{eq:parametric_nlp} with $x = x^{k+1}_{\mathrm{pred}}$.
Then, in \ref{step:feedback}, a linear system or QP is solved to obtain $z^{k+1}$~\cite{Zavala2009,Jaeschke2014}.
If an advanced problem with a perfect prediction $x^{k+1} = x^{k+1}_{\mathrm{pred}}$ is solved to local optimality, there is no numerical error in the feedback, i.e., $z^{k+1} = \bar{z}(x^{k+1})$.
We denote AS-RTI with level X iteration as \textit{AS-RTI-X}.

\vspace{-0.5cm}
\subsection{Multi-Level Iterations}\label{sec:mli}
\vspace{-0.04cm}
In step~\ref{step:preparation}, we use some Multi-Level Iteration (MLI) variant to compute $z^{k}_{\lin}$.
All MLI levels start with a reference point $\hat{z}^k = (\hat{w}^k,\hat{\lambda}^k,\hat{\mu}^k)$, which can be e.g., the linearization point for the previous feedback phase in the AS-RTI context.
The different MLI levels recompute different values of \eqref{eq:qp_sqp} and use evaluations at $\hat{z}^k$ for the others.

\textit{Level D iterations.}
Level D iterations are essentially full SQP iterations as described in Section \ref{sec:sqp}.
All functions and derivative derivations are evaluated exactly.

\textit{Level C iterations.}
In level C iterations, all matrices in QP~\eqref{eq:qp_sqp} are fixed:
 $\hat{A}^k \approx \nabla_{ww}^2 \mathcal{L}(\hat{z}^k)$, $\hat{G}^k = \nabla_w g(\hat{w}^k)^\top$, $\hat{H}^k = \nabla_w g(\hat{w}^k)^\top$.
Only the vectors, i.e., $g^{k,j}$, $h^{k,j}$ and $a^{k,j}$, are updated.
Note that the SQP subproblem \eqref{eq:qp_sqp} only uses the objective gradient instead of the Lagrangian for $a^{k,j}$, since primal variables are updated in a delta and the duals in an absolute fashion as in \cite[Sec. 18.1]{Nocedal2006}.
To modify \eqref{eq:qp_sqp} to take the latest multipliers into account in the Lagrange gradient but fix the linearization of the constraints, $a^{k,j}$ needs to be updated for $j\geq0$ as
	$a^{k,j} = \nabla\mathcal{L}(w^{k,j},\lambda^{k,j},\mu^{k,j}) + \lambda^{k,j} \hat{G}^k + \mu^{k,j} \hat{H}^k.$
For a $\hat{z}^k$ close enough to a solution, these iterations converge linearly to a local optimum of \eqref{eq:parametric_nlp}~\cite{Bock2007}.

\textit{Level B iterations.}
In level B iterations, also called zero-order iterations~\cite{Zanelli2021c}, only functions and no derivatives are evaluated to set up QP~\eqref{eq:qp_sqp}.
More specifically, only $g^{k,j}, h^{k,j}$ are obtained by evaluations.
The objective gradient is approximated by $a^{k,j} = \nabla_w f(\hat{w}^k) + \hat{A}({w}^{k,j}_\B-\hat{w}^k)$.
Here, $w_\B^{k,j+1}$ are the iterates computed by $w_{\B}^{k,j+1} = w_{\B}^{k,j}+\Delta w^{k,j}, w_{\B}^{k,0} = \hat{w}^k$, $\lambda_{\B}^{k,j+1}= \lambda^{k,j}_{\mathrm{QP}}$ and $\mu_{\B}^{k,j+1}= \mu^{k,j}_{\mathrm{QP}}$, where $(\Delta w^{k,j}, \lambda^{k,j}_{\mathrm{QP}}, \mu^{k,j}_{\mathrm{QP}})$ is the solution of QP~\eqref{eq:qp_sqp} with the partially updated data as just described.
We call the generated sequence $\{z^{k,j}_\B\}_{j\geq0}$ for a fixed $k$ the \textit{level B iterates}.

\begin{proposition}(Adapted from \cite[Theorem 1.4]{Bock2007}) \label{prop:level_B_moritz}
	If for a fixed parameter $x^k$, the level B iterates $\{z^{k,j}_{\B}\}$ converge to a limit point $\bar{z}_{\mathrm{B}}(x^k) = (\bar{w}_{\B}(x^k), \bar{\lambda}_{\B}(x^k), \bar{\mu}_{\B}(x^k))$, short $\bar{z}_\B^k$, then $\bar{z}_{\B}^k$ is a primal-dual solution of the parametric NLP
	\vspace{-0.1cm}
\begin{mini!}
	{\substack{w}}
	{f(w)+w^\top \beta^k}
	{\label{eq:parametric_nlp_B}}
	{}
	\addConstraint{0}{= g(w) + Mx^k\label{eq:parametric_nlp_eq_B}}
	\addConstraint{0}{\leq h(w), \label{eq:parametric_nlp_ineq_B}}
\end{mini!}
\vspace{-.6cm}
\begin{align}
	\text{with} \;
	\beta^k &= \nabla f(\hat{w}^k)+\hat{A}^k(\bar{w}_{\B}^k-\hat{w}^k)-\nabla f(\bar{w}_{\B}^k)\nonumber \\
	&+(\nabla g(\bar{w}_{\B}^k)\shortminus (\hat{G}^k)^\top) \bar{\lambda}_{\B}^k
	+(\nabla h(\bar{w}_{\B}^k)\shortminus (\hat{H}^k)^\top )\bar{\mu}_{\B}^k. \nonumber
\end{align}
\end{proposition}
Therefore, the level B iterations converge to a solution of an NLP that is parametric in $x^k$ and $\beta^k$.
The local minimizer of \eqref{eq:parametric_nlp_B} is feasible for \eqref{eq:parametric_nlp} but not optimal, as the objective is altered by $w^\top \beta^k$.

\textit{Level A iterations.}
Compared to level B iterations, a level A iteration does not evaluate the constraint residuals and only updates them with respect to the parameter, i.e., $g^{k} = g(\hat{w}^k)+x^k, h^k = h(\hat{w}^k)$.
The new solution approximation is obtained by a single QP solve, with $w_{\A}^k = \hat{w}^k + \Delta w$.
It is important to note that level A iterations only generate a new value $w_{\A}^k$ for every new parameter value $x^k$, instead of a sequence as in the other levels.
The reason is that the QP is a piecewise linear approximation of the solution map $\bar{z}(x)$, which is evaluated at given parameters $x^k$.

\vspace{-.2cm}
\section{Improved error estimates for the AS-RTI}\label{sec:theory}
This section recalls some results for the convergence of predictor-corrector methods and derives novel error bounds for different variations of the AS-RTI method.
\vspace{-0.3cm}
\subsection{Error bounds for the feedback phase}
For ease of exposition, let us consider only equality-constrained problems, or assume a fixed active set.
Later, in Remark~\ref{remark:inequalities}, we comment on what needs to be changed to generalize the result to inequality-constrained problems.
Regard the parametric optimization problem
\begin{align}\label{eq:parametric_nlp_simple}
	\min_{{w}} f(w) \quad \mathrm{s.t.} \quad 0 = g(w)+Mx.
\end{align}

The Karush–Kuhn–Tucker (KKT) conditions of~\eqref{eq:parametric_nlp_simple} can be stated as the nonlinear root-finding problem
\begin{align}\label{eq:kkt_conditions}
	F(z,x) = \begin{bmatrix}
		\nabla f(w) - \nabla g(w) \lambda\\
		g(w) + Mx
	\end{bmatrix}= 0,
\end{align}
with $z = (w,\lambda)\in \R^{n_z}$, $n_z = n_w+n_g$ and the function $F:\R^{n_z}\times \R^{n_x} \to \R^{n_z}$ which is at least once continuously differentiable.
Its zeros are denoted by $\bar{z}(x)= \bar{z}$.

We make a regularity assumption on the solutions of the parametric NLP \eqref{eq:parametric_nlp_simple}.
The Linear Independence Constraint Qualification (LICQ) is said to hold at a point $w$ if the vectors $\nabla g_i(w), i =1,\ldots,n_g$ are linearly independent.
The Second order Sufficient Conditions (SOSC) is said to hold at a KKT point $\bar{z}^k$  if $Z^\top \nabla_{ww}^2 \mathcal{L}(\bar{z}^k)Z \succ 0$, where $Z\in \R^{n_w\times (n_w-n_g)}$ is a basis for the null space of~$\nabla g(\bar{z}^k)^\top$.

\begin{assumption}(LICQ, SOSC)\label{ass:licq_sosc}
For all parameters $x \in X \subseteq \R^{n_x}$,  all local minimizers $\bar{z}(x)$ of \eqref{eq:parametric_nlp} satisfy the LICQ and SOSC conditions.
\end{assumption}

The Jacobian of $F(z,x)$, which due to linearity in $x$ only depends on $z$, is denoted by $J(z) \coloneqq \frac{\partial (F(z,x))}{\partial z}$.
Recall that $z^{k+1} \approx \bar{z}(x^{k+1})$ and $z^{k} \approx \bar{z}(x^{k})$.
A predictor-corrector step for the parametric root-finding problem \eqref{eq:kkt_conditions} can be written
\begin{align}\label{eq:predictor_corrector}
			z^{k+1} & =z^{k} - J(z^{k})^{-1} F(z^k,x^{k+1}).
\end{align}
Note the dependence on the next parameter $x^{k+1}$.
If a fixed $x^k$ is used instead, equation~\eqref{eq:predictor_corrector} reduces to an exact Newton-step, a corrector step.
If $z^{k} = \bar{z}(x^k)$, due to the linearity of $F(z,x)$ in $x$, equation~\eqref{eq:predictor_corrector} reduces to a predictor step: $z^{k+1} = z^k  + \frac{\partial \bar{z}(x^{k})}{\partial x}(x^{k+1}-x^k)$, see \cite[Sec. 1.2]{TranDinh2012b} for details.

\color{purple}
\color{black}

We proceed by restating some results that we use to study the AS-RTI error, starting with the convergence of Newton's method, cf. \cite[Sec. 8.3.3.]{Rawlings2017}.
Suppose that the parameter $x^k$ is fixed, then \eqref{eq:predictor_corrector} is simply a standard exact Newton-step
$z^{k,j+1} = z^{k,j} - J(z^{k,j})^{-1}F(z^{k,j},x^{k})$.
One may also use a Jacobian approximation $M^{k,j}\approx J(z^{k,j})$.
\vspace{-.1cm}
\begin{assumption}\label{ass:omega_kappa}
	($\omega$ and $\kappa$ conditions)
	There exist $\omega < \infty$ and $\kappa < 1$ such that, for any fixed $x^k$, given iterate $z^{k,j}$ and Jacobian approximation $M^{k,j}$, the following holds:
	\begin{enumerate}[(a)]
		\item $ \|(M^{k,j})^{-1}( J(z^{k,j}) - J(z))\| \leq \omega \|{z^{k,j} - z}\|,\ \forall z, $
		\item $ \|(M^{k,j})^{-1}( J(z^{k,j}) - M^{k,j})\| \leq \kappa$.
	\end{enumerate}
\end{assumption}
Here, (a) is a rescaled Lipschitz condition on the Jacobian, and (b) measures the Jacobian approximation error.
For the exact Jacobians $(M^{k,j}= J(z^{k,j}))$ it holds that ${\kappa = 0}$.
\vspace{-.1cm}
\begin{assumption}\label{ass:initialization}(Initialization)
	A starting point $z^{k,0}$  of a sequence generated by a Newton-type method satisfies
	\begin{align}\label{eq:inital_pint}
		\| \bar{z}^k-z^{k,0} \| <r_z =  \tfrac{2(1-\kappa)}{\omega}.
	\end{align}
\end{assumption}
Assumption~\ref{ass:initialization} tells how close to a solution one must initialize so that Newton's method converges with full steps.
Both Ass.~\ref{ass:omega_kappa} and~\ref{ass:initialization} are standard for local Newton-type convergence analysis to state the following theorem~\cite{Diehl2001, Bock2007, Zanelli2021c, TranDinh2012b}.

\begin{theorem}(Newton-type convergence, cf. \cite[Thm. 8.7]{Rawlings2017}) \label{th:newton_convergnece}
	Regard the continuously differentiable function $F$ in \eqref{eq:kkt_conditions} with a fixed $x^k$, and a solution $\bar{z}(x^k)$ (short $\bar{z}^k$) with $F(\bar{z}^k, x^k)=0$.
	If Ass. \ref{ass:omega_kappa} holds, the sequence $\{z^{k,j}\}$ generated by $z^{k,j+1} = z^{k,j} - (M^{k,j})^{-1}F(z^{k,j},x^k)$ satisfies the inequality
	\vspace{-.1cm}
	\begin{align} \label{eq:newton_convergence}
		\| \bar{z}^k- z^{k,j+1}\| &\leq \left(\kappa + \tfrac{\omega}{2}	 \| \bar{z}^k- z^{k,j}\|  \right)\| \bar{z}^k- z^{k,j}\|.
	\end{align}
	Moreover, if Assumption \ref{ass:initialization} holds,
	then the sequence $\{z^{k,j}\}$ converges to $\bar{z}^k$ as $j\to\infty$.
\end{theorem}
Define $\alpha^k \coloneqq \kappa + \frac{\omega}{2} \| \bar{z}^k- z^{k,0}\|$, with $\alpha^k<1$ due to Ass.~\ref{ass:initialization}.
Applying \eqref{eq:newton_convergence} recursively, results in the useful inequality
\begin{align} \label{eq:newton_convergence_0}
	\| \bar{z}^k- z^{k,j+1}\| &\leq (\alpha^k)^{j+1} \| \bar{z}- z^{k,0}\|.
\end{align}

Next, we look at the solution map $\bar{z}(x)$ for different parameters $x$.
Assumption \ref{ass:licq_sosc} implies that $J(\bar{z}(x))$ is invertible for all $x\in X$, cf. \cite[Lem. 16.1]{Nocedal2006}.
Thus, applying the implicit function theorem \cite[Thm. A.1]{Nocedal2006} to \eqref{eq:kkt_conditions} entails that $\bar{z}(x)$ is locally unique and Lipschitz continuous
\begin{align}\label{eq:lipschitz_z}
	\| \bar{z}(x^{k+1}) - \bar{z}(x^{k}) \| \leq \sigma \| x^{k+1}-x^{k}\|.
\end{align}
By applying a single predictor-corrector step \eqref{eq:predictor_corrector} for the parameters $x^{k+1}$ and $x^{k}_{\mathrm{lin}}$, the corresponding solution map approximations $z^{k+1}\!\approx\!\bar{z}(x^{k+1})$ and linearization point $z^{k}_{\mathrm{lin}} \approx \bar{z}(x^{k}_{\mathrm{lin}})$ can be related by combining \eqref{eq:newton_convergence} and \eqref{eq:lipschitz_z} as follows.
\begin{theorem}(Adapted from \cite[Th. 3.5]{TranDinh2012b}, \cite[Lem. 3.1.5]{Zanelli2021c})\label{th:predictor_corrector}
	Regard the continuously differentiable function $F$ in \eqref{eq:kkt_conditions}, and solution $\bar{z}(x)$ with $F(\bar{z}(x),x)=0$.
	Let Assumptions~\ref{ass:licq_sosc} and \ref{ass:omega_kappa} hold.
	Then, the iterates $z^{k+1}$ and $\bar{z}^{k}_{\mathrm{lin}}=z^k$ generated by~\eqref{eq:predictor_corrector} for the parameters $x^{k+1}$ and $x^{k}_{\mathrm{lin}}=x^k$ satisfy
\vspace{-.1cm}
	\begin{align} \label{eq:predictor_corrector_bound}
			\| &\bar{z}^{k+1}\!-\!{z}^{k+1}\| \leq
			(\kappa \sigma + \tfrac{\omega \sigma^2}{2}\|x^{k+1} - x^{k}_{\mathrm{lin}}\| ) \|x^{k+1} \!-\! x^{k}_{\mathrm{lin}}\|\nonumber\\
			&+\!(\kappa\!+\!\omega \sigma \|x^{k+1} \!\shortminus \! x^{k}_{\mathrm{lin}}\| \!+\! \tfrac{\omega}{2}\|\bar{z}^{k}_{\mathrm{lin}} \!\shortminus\! {z}^{k}_{\mathrm{lin}}\|) \|\bar{z}^{k}_{\mathrm{lin}} \!\shortminus \! {z}^{k}_{\mathrm{lin}}\|.
	\end{align}
\vspace{-.1cm}
\end{theorem}
Note that index $j$ is omitted since a single Newton-type step is computed for every parameter.
Furthermore, if for some fixed $x^0$ the point $z^0 \approx \bar{z}^0$ satisfies Ass.~\ref{ass:initialization}, there exists $r_x$ (depending on $\sigma, \kappa$ and $\omega$) such that for $\|x^{k+1} \shortminus x^k\| < r_x$, the sequence $\{z^k\}$ generated by \eqref{eq:predictor_corrector} remains bounded with $\|z^k - \bar{z}(x^k)\| < r_z$, cf. \cite[Corollary 3.6.]{TranDinh2012b}.

\vspace{-0.3cm}
\subsection{Error bounds for the preparation phase} \label{sec:theory_prep}
Given a linearization point ${z}^k_{\mathrm{lin}}$ and $x^k_{\mathrm{lin}}$, in the feedback phase, i.e., step \ref{step:feedback} of AS-RTI, a QP is solved to obtain an approximation $z^{k+1} \approx \bar{z}(x^{k+1})$.
In the absence of inequalities, instead of a QP, the linear system \eqref{eq:predictor_corrector} is solved, and we have the error bound in \eqref{eq:predictor_corrector_bound}.
We see from \eqref{eq:predictor_corrector_bound} that the accuracy of the new output $z^{k+1}$ improves with a smaller difference between the parameters via $\|x^{k+1} - x^{k}_{\mathrm{lin}}\|$ and higher accuracy of the linearization point via $\|\bar{z}^{k}_{\mathrm{lin}}-{z}^{k}_{\mathrm{lin}}\|$.
The goal of AS-RTI is to reduce $\|\bar{z}^{k}_{\mathrm{lin}}-{z}^{k}_{\mathrm{lin}}\|$ by solving an advanced problem with $x^{k}_{\mathrm{lin}} = \xpred \approx x^{k+1}$ performing iterations with some MLI variant, which results in tighter bounds for $\|\bar{z}^{k+1}\!-\!{z}^{k+1}\| $ in~\eqref{eq:predictor_corrector_bound}.
We proceed by quantifying the error for each variant.

\textit{Level D iterations.}
Here, we use Newton-type steps in~\eqref{eq:predictor_corrector} with $M^{k,j}_{\mathrm{D}} \approx J(z^{k,j})$ and assume that Ass. \ref{ass:omega_kappa} holds for some $\kappa_\D$ and $\omega_\D$.
Starting with the previous output $z^k$ and predicted parameter $\xpred$, and assuming that $\| z^k - \bar{z}(\xpred)\|< \frac{2(1-\kappa_{\D})}{\omega_\D}$, a modification of Ass. \ref{ass:initialization}, we carry out $N_\D$ Newton-type iterations
Theorem~\ref{th:newton_convergnece} and \eqref{eq:newton_convergence_0} yield
\begin{align}
	\label{eq:contract_D}
	\| z_\lin^{k} - \bar{z}(\xpred) \| \leq (\alpha_\D^k)^{N_\D} \| z^k - \bar{z}(\xpred)\| .
\end{align}
with $\alpha_\D^k = \kappa_\D + \frac{\omega_\D}{2}\| z^k - \bar{z}(\xpred)\| < 1$.
In the limiting case where $\xpred = x^{k+1}$ and $j \to \infty$, we see that in the right-hand side of \eqref{eq:predictor_corrector_bound} becomes zero, i.e. $z^{k+1} = \bar{z}^{k+1}$.

\textit{Level C iterations.}
In level C, we essentially proceed as in level D, except that we use a constant Jacobian approximation $M^{k,j} = M^k_\C = \begin{bsmallmatrix}
\hat{A}^k & (\hat{G}^k)^\top \\
\hat{G}^k & 0\end{bsmallmatrix}$ computed at the reference point $\hat{z}^k=z^{k-1}$.
We assume that that Ass.~\ref{ass:omega_kappa} holds for some $\kappa_{\C}$ and $\omega_\C$.
Similar to level D, starting with the previous output, $z^k$ and $\xpred$, and assuming that $\| z^k - \bar{z}(\xpred)\|< r_{z,\C}  = \frac{2(1-\kappa_\C)}{\omega_\C}$, a modification of Ass.~\ref{ass:initialization}, we perform $N_\C$ Newton-type iterations
to obtain $z^k_{\lin}$.
Applying Thm~\ref{th:newton_convergnece} and \eqref{eq:newton_convergence_0} yields
\begin{align}
	\label{eq:contract_C}
	\| z_\lin^{k} - \bar{z}(\xpred) \| \leq (\alpha_\C^k)^{N_\C} \| z^k - \bar{z}(\xpred)\| .
\end{align}
with $\alpha_\C^k = \kappa_\C + \frac{\omega_\C}{2}\| z^k - \bar{z}(\xpred)\| < 1$.
Since $M_\C^k$ is fixed, it is a less accurate but computationally cheaper Jacobian approximation and usually $\kappa_\D \ll \kappa_\C$ holds.
For similar $\omega_\D$ and $\omega_\C$, it follows that $r_{z,\C} < r_{z,\D}$.
In other words, the previous output $z^k$ must be closer to the solution $\bar{z}(\xpred)$ for level C than for D to achieve contraction of the iterates.

\textit{Level B iterations.}
Following Proposition~\ref{prop:level_B_moritz}, Level B iterations for NLPs without inequality constraints converge to a solution of
\begin{align}\label{eq:parametric_nlp_simple_B}
	\min_{{w}} f(w) +  w^\top \beta^k \quad \mathrm{s.t.} \quad 0 = g(w)+Mx^k,
\end{align}
with $\beta^k = \nabla f(\hat{w}^k)+\hat{A}^k(\bar{w}_{\B}^k-\hat{w}^k)-\nabla f(\bar{w}_{\B}^k)+(\nabla g(\bar{w}_{\B})^k-(\hat{G}^k)^\top) \bar{\lambda}_{\B}^k$.
The KKT conditions of \eqref{eq:parametric_nlp_simple_B} read
\begin{align}\label{eq:kkt_conditions_B}
	F_\B(z,\beta^k,x^k)
	= F(z,x^k) + [(\beta^k)\transp, 0]\transp.
\end{align}
A sequence of iterates is generated via $z^{k,j+1} = z^{k,j} - (M^k_\B)^{-1} F_\B(z^j,\beta,\xpred)$, where the parameters $\beta^k$ and $\xpred$ are fixed, and $M^k_\B = M^k_\C$.
We have the following estimate.
\begin{proposition}\label{prop:level_B_armin}
	Assume that LICQ and SOSC hold for~\eqref{eq:parametric_nlp_simple_B} at all $\bar{z}_B(x,\beta)$, for $x\in X$ and $\beta\in \R^{n_w}$.
	Suppose that Ass.~\ref{ass:omega_kappa} holds for $z^{k,j+1} = z^{k,j} - (M_\B^k)^{-1} F_\B(z^j,\beta^k,\xpred)$ with constants $\kappa_\B$ and $\omega_B$, and that $\| z^k - \bar{z}^k_\B\| < \frac{2(1-\kappa_\B)}{\omega_\B}$.
	Then the sequence of iterates $\{z_\B^j\}_{j=1,\ldots, N_\B}$ fulfills
	\begin{align*}
		\| \bar{z}(\xpred)\shortminus z_\lin^k \| \leq \sigma_\B \| \beta^k \| + (\alpha^k_\B)^{N_\B} \| \bar{z}_\B(\xpred,\beta^k)\!  \shortminus  \! z^k\|,
	\end{align*}
	where $\sigma_B$ is the Lipschitz constant of $\bar{z}_\B(x,\beta)$ and  $\alpha_\B^k = \kappa_\B + \frac{\omega_B}{2}\| z^k - \bar{z}_B(\xpred,\beta^k)\| < 1$.
\end{proposition}
\textit{Proof.}
Adding and subtracting $\bar{z}_\B(x,\beta^k)$ in the left term of the next equation, and using the triangle inequality we get
\begin{align*}
	\| \bar{z}(x)-z_\lin^k \| \leq \| \bar{z}(x)-\bar{z}_\B(x,\beta^k) \|  + \|\bar{z}_\B(x,\beta^k)-z_\lin^k \|.
\end{align*}
For the first term on the right, we note that $\bar{z}_\B(x,0) = \bar{z}(x)$.
Since LICQ, SOSC hold for \eqref{eq:parametric_nlp_simple_B} by applying the implicit function theorem to \eqref{eq:kkt_conditions_B}, we have that
$\| \bar{z}(x)-\bar{z}_\B(x,\beta) \| \leq \sigma_\B \|\beta\|$.
For the second term on the left, we apply~\eqref{eq:newton_convergence_0} with $N_\B$ iterations, and obtain that $\|\bar{z}_\B(x,\beta)-z_\lin^k \| \leq (\alpha_\B)^{N_\B} \| \bar{z}_\B(\xpred,\beta)  - z^k\|$.
 By using these two terms and $x = \xpred$ the result of this proposition is obtained.
 \qed

Even if level B iterations are fully converged $j\to \infty$, the error $\| \bar{z}(\xpred)-z_\lin^k \| \leq \sigma_\B \| \beta \|$ remains, as they converge to feasible but suboptimal points of \eqref{eq:parametric_nlp}, cf. Prop.~\ref{prop:level_B_moritz}.
Analogously to C and D, the previous output $z^k$ has to be close enough to $\bar{z}_B(\xpred)$ such that the iterates contract.

\textit{Level A iterations.} %
Applying Thm~\ref{th:predictor_corrector} with $x^{k}_\lin\!=\!\xpred$ yields an error bound for the feedback phase of AS-RTI-A.
\begin{proposition}\label{prop:level_A}
Let the assumptions of Theorem~\ref{th:predictor_corrector} hold.
When ${z}^{k}_{\mathrm{lin}}$ is obtained via a Level A iteration ${z}^{k}_{\mathrm{lin}} = z^{k\shortminus 1} \shortminus (M^{k\shortminus 1})^{\shortminus 1} F(z^{k\shortminus 1},\xpred)$, the feedback error satisfies
\begin{align}
	\| &\bar{z}^{k+1}\!\shortminus\!{z}^{k+1}\| \leq
	(\kappa \sigma + \tfrac{\omega \sigma^2}{2}\|x^{k+1} \shortminus \xpred \| ) \|x^{k+1} \!\shortminus\! \xpred \|\nonumber\\
	&+\!(\kappa\!+\!\omega \sigma \|x^{k+1} \!\shortminus \! \xpred \| \!+\! \tfrac{\omega}{2}\|\bar{z}^{k}_{\mathrm{lin}} \!\shortminus\! {z}^{k}_{\mathrm{lin}}\|) \|\bar{z}^{k}_{\mathrm{lin}} \!\shortminus \! {z}^{k}_{\mathrm{lin}}\|.
\end{align}
\end{proposition}

Compared to RTI, where $x_\lin^k\!=\!x^k$ and $z^k_\lin\!=\!z^k$, AS-RTI-A can significantly reduce the error if $\xpred \approx x^{k+1}$ and if the predictor provides a good approximation~${z}^{k}_{\mathrm{lin}}\!\approx\!\bar{z}(\xpred)\!=\! \bar{z}^{k}_{\mathrm{lin}} $.
\begin{remark}
The same linearization at $z^{k\shortminus 1}$ is reused to compute both an approximation for $\xpred$ and $x^k$.
It is important to start from $z^{k\shortminus 1}$ for computing $z^k_\lin$ instead of $z^k$ to avoid taking the same corrector step twice, cf. \eqref{eq:predictor_corrector}.
\end{remark}

\begin{remark}[Extension to inequality constraints]
\label{remark:inequalities}
All results are derived from two inequalities: convergence of Newton's method in Theorem~\ref{th:newton_convergnece}, and Lipschitz continuity of $\bar{z}(x)$ in~\eqref{eq:lipschitz_z}.
The KKT conditions of inequality constrained NLP~\eqref{eq:parametric_nlp} can be written as a generalized equation solved by the Newton-Josephy method, which is equivalent to SQP.
Using Robinson's strong regularity, which is implied by strong SOSC and LICQ, Theorem~\ref{th:newton_convergnece} can be generalized via \cite[Thm. 3.5]{TranDinh2012b} and Eq.~\eqref{eq:lipschitz_z} via \cite[Lem. 3.3]{TranDinh2012b}.
\end{remark}

\vspace{-.35cm}
\section{Implementation in \acados}\label{sec:implementation}
Next, we discuss some practical aspects for an efficient implementation of AS-RTI within the \acados{} software, which has been developed as part of this work.

\paragraph{Condensing and QP solution with two phases}
The \acados{} software offers a variety of QP solvers.
HPIPM offers efficient methods to transform OCP-structured QPs into dense ones or ones with a shorter horizon by full and partial condensing \cite{Frison2016}.
Since in the preparation phase all matrices of the QP are readily available, most of the condensing operations can be performed in that phase.
Moreover, it is common to eliminate the initial state variable from the QP.
An efficient split of operations is realized by implementing functions that assume that only matrices of QP \eqref{eq:qp_sqp} are known and a second one that completes the computations once the vector quantities are known.
This split functionality is utilized when implementing the level A, B, and C iterations.

\paragraph{Advancing}
There are two main strategies to set up the advanced problem,~\ref{step:advance}.
1) Simulate with $\phi_0(\cdot)$ internally at the current SQP iterate.
2) Simulate the system externally.
Both strategies work in general for nonuniform discretization grids in \eqref{eq:acados_OCP}, which have been shown to be superior with respect to uniform ones in \cite{Frey2023b}.
Option 1) assumes that $\phi_0$ models the systems evolution over the sampling time.
In contrast, 2) does not require this assumption and allows using a higher fidelity model than used in the OCP.

\vspace{-0.25cm}
\section{Numerical experiments}
\label{sec:experiments}

\setlength{\tabcolsep}{5pt}
\begin{table}
	\vspace{0.2cm}
	\centering
	\caption{Timings, relative suboptimality, stationarity residual and constraint violation for different controllers. \label{tab:asrti_experiment}}
	\vspace{-.2cm}{\footnotesize
{\footnotesize
	\begin{tabular}{lrrrrrr}
	\toprule
	& max. time & max. time & relative & mean & mean \\
	& prepare & feedback & subopt. & $10^3 \norm{g} $ & $\norm{\nabla_w\mathcal{L}}$ \\
	algorithm & [ms] & [ms] &  [\%] &  &  \\ \midrule
	SQP-100&{0.00} &5.518 &{0.05} &0.00 &0.00\\
	SQP-2&{0.00} &0.264 &{0.25} &12.99 &7.75\\
	AS-RTI-D-2&0.35 &{0.020} &{0.04} &0.74 &1.28\\
	AS-RTI-D-1&0.23 &{0.020} &{0.25} &8.26 &5.38\\
	AS-RTI-C-2&0.32 &{0.020} &{0.04} &5.34 &2.20\\
	AS-RTI-C-1&0.23 &{0.021} &{0.24} &10.70 &5.85\\
	AS-RTI-B-2&0.28 &{0.021} &{0.57} &1.28 &8.20\\
	AS-RTI-B-1&0.19 &{0.021} &{0.54} &8.77 &7.98\\
	AS-RTI-A&0.13 &{0.019} &{0.54} &12.41 &6.63\\
	RTI&0.11 &{0.022} &3.55 &14.06 &8.26\\
	\bottomrule
	\end{tabular}
	}
	\vspace{-.3cm}}
\end{table}

\begin{figure}
	\includegraphics[width=\columnwidth]{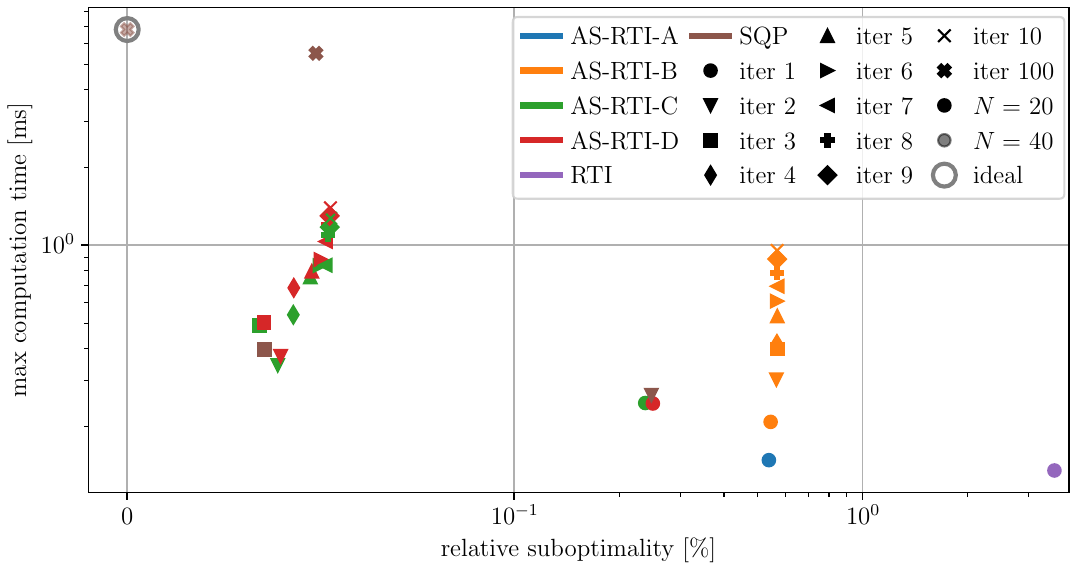}
	\vspace{-.7cm}
	\caption{Pareto plot: timings vs. relative suboptimality.
	\label{fig:pareto_pendulum_subopt}
	}
	\vspace{-.6cm}
\end{figure}

This section compares different MPC controllers in an open-source numerical simulation study \cite{as_rti_acados} using \href{https://github.com/acados/acados/releases/tag/v0.3.2}{\acados{} v0.3.2} \cite{acados_releases}
via Python on a Laptop with an Intel i5-8365U CPU, 16 GB of RAM running Ubuntu 22.04.

\paragraph{Inverted pendulum on cart test problem}
The differential state of the model is $x = [p, \theta, s, \omega]^\top$ with cart position $p$, cart velocity $s$, angle of the pendulum $\theta$ and angular velocity $\omega$.
The control input $u$ is a force acting on the cart in the horizontal plane.
The ODE, describing the system dynamics can be found e.g. in~\cite{Verschueren2021}.
In our OCP formulation, $u$ is constrained to be in $[-40, 40]$.
The goal is to stabilize the system in the unstable upright position driving all states to zero.
We formulate the linear least squares cost $l(x, u) = x^\top Q x + u^\top R u $ with cost weights are $Q = \mathrm{diag}(100, 10^3, 0.01, 0.01) $, $ R = 0.2 $.
The terminal cost term is set to $E(x) = x^\top P x$,
where $P$ is obtained as a solution of the discrete algebraic Riccati equation with cost and dynamics linearized at the steady state.

\paragraph{Scenario}
The system is simulated for four seconds at a sampling time of $\Delta t = 0.05 \mathrm{s}$.
We simulate 20 different scenarios, in each of which, the system starts at an upward position with a random initial value for $p$ and is disturbed at two time instances, at $0\mathrm{s}$ and $2\mathrm{s}$ by overwriting the control action with random value in $[-100, 100]$.
The OCP is formulated with a time horizon of $2\mathrm{s}$ divided into $N=20$ shooting intervals, the first is of length $\Delta t$ and the remaining uniformly split the rest of the time horizon.
The dynamics are discretized using one step of an implicit Radau IIA method of order three with three Newton iterations on each shooting interval, respectively of order seven with 20 Newton iterations for the simulation step of the plant.

\paragraph{Controllers}
We apply a variety of controllers and report their performance in Table~\ref{tab:asrti_experiment}.
All controllers use the full condensing functionality from \texttt{HPIPM}~\cite{Frison2016} and the active-set solver \texttt{DAQP}~\cite{Arnstrom2022}.
The solvers labeled \textit{SQP~n} apply $n$ SQP iterations.
Different variants of AS-RTI controllers are labeled \textit{AS-RTI-X-$n$} performing $n$ level X iterations on the advanced problem in each preparation phase.
Additionally, we compare with a plain \textit{RTI} controller.

\begin{figure}
	\vspace{.2cm}
	\includegraphics[width=\columnwidth]{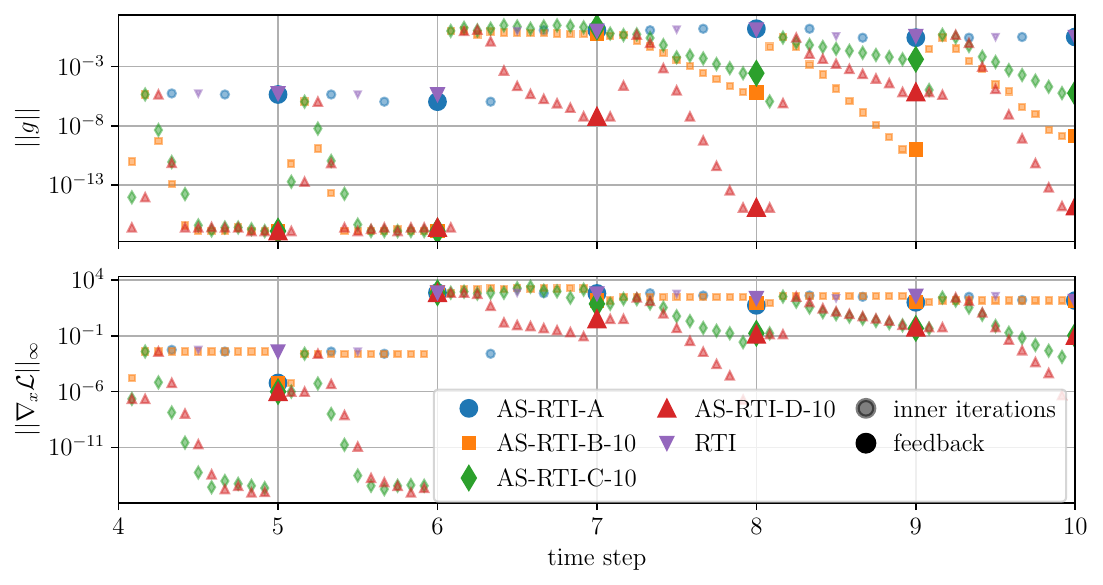}
	\vspace{-.8cm}
	\caption{Primal and dual residuals of inner iterations and the applied feedback steps for RTI and different AS-RTI variants.
	\vspace{-1.3cm}
	\label{fig:infeasibility_pendulum_as_rti}
	}
\end{figure}

\paragraph{Evaluation}
The Pareto plot in Fig.~\ref{fig:pareto_pendulum_subopt} compares the controller variants in terms of maximum computation time and relative suboptimality.
The latter is evaluated by computing the closed-loop cost and comparing with a controller that uses a finer uniform discretization grid with $N=40$ and fully converged SQP, which is marked as \textit{ideal} in Fig.~\ref{fig:pareto_pendulum_subopt}.
In Table~\ref{tab:asrti_experiment}, we additionally report the maximum timings for preparation and feedback phase over all simulations.
Additionally, the mean values of the constraint violation and the Lagrange gradient over all simulation steps are listed.
The shooting gaps show satisfaction of nonlinear constraints, quantified as $\norm{g}$.
The example only contains linear inequality constraints, which are always satisfied.
Thus, it allows one to compare closed-loop suboptimality and constraint violation simultaneously.
Figure~\ref{fig:infeasibility_pendulum_as_rti} shows how primal and dual infeasibility evolve over all AS-RTI iterations over a few time steps after applying a large disturbance to the plant.

\paragraph{Discussion}
Firstly, we can see from Table~\ref{tab:asrti_experiment} that the timings of the feedback step for all AS-RTI variations are consistent and a multiple lower compared to their preparation phase, which is enabled by the split condensing described in Sec.~\ref{sec:implementation}.
In our closed-loop simulation, the delay from the computation of the feedback phase is neglected, which would impact all real-time controllers similarly, but would drastically degrade the performance of the SQP algorithms included for reference.
The classic RTI algorithm has a rather high suboptimality and is the fastest real-time algorithm.

AS-RTI-A, with a single precondensed QP solve, is able to greatly improve on the classic RTI algorithm with much better performance and a marginally increased computational load, cf. Fig.~\ref{fig:pareto_pendulum_subopt}.
In Figure~\ref{fig:infeasibility_pendulum_as_rti}, one can observe that level~B iterations converge to a feasible linearization point.
This iterate is associated with a fixed suboptimality which is visible in the second subplot and consistent Proposition~\ref{prop:level_B_armin}.

In contrast to the level B iterations, the level C iterations converge to a feasible locally optimal point.
The same holds for level D iterations.
Since iterations of level D are more accurate than level C, it follows that $\kappa_\D < \kappa_\C$, which implies a faster error reduction in \eqref{eq:contract_D} compared to \eqref{eq:contract_C}.
This can be observed in Figure~\ref{fig:infeasibility_pendulum_as_rti} which shows faster convergence for D, i.e., the steeper slope in both primal and dual infeasibility.

Overall, we see in Figure~\ref{fig:infeasibility_pendulum_as_rti} that the feedback iterates of all AS-RTI variants are more accurate than the one of plain RTI.
This is due to the tighter bound in \eqref{eq:predictor_corrector_bound}, more precisely, the reduced error in the linearization point $\|\bar{z}^{k}_{\mathrm{lin}} -  {z}^{k}_{\mathrm{lin}}\| $, and the smaller parameter difference $ \|x^{k+1} - x^{k+1}_{\mathrm{pred}}\| $ for AS-RTI compared to $ \|x^{k+1} - x^{k}\| $ for plain RTI.

\section{Conclusion}\label{sec:conclusion}
This paper streamlines and extends the existing analysis of AS-RTI with Multi-Level Iterations (MLI) of all for levels.
It is shown that if the current solution is sufficiently close to the next solution, the numerical error can be reduced by a few computationally cheap MLI iterations.
Furthermore, an efficient implementation in the open-source package~\acados{} is presented making the method widely available for real-world applications on embedded hardware.
Numerical examples confirm the theory and demonstrate how to assemble efficient variants of the AS-RTI method.
In particular, AS-RTI-A can significantly improve control performance over standard RTI at only 20\% additional cost per sampling time, by computing a single additional QP solution to the advanced problem.

\bibliography{syscop}

\end{document}

%% file: oj_commands.tex
\usepackage{mathtools}

\definecolor{wheat}{rgb}{0.96,0.87,0.70}
\newcommand{\transp}{^\top}

\newcommand{\R}{\mathbb{R}}

\newcommand{\acados}{\texttt{acados}}

\DeclarePairedDelimiter\abs{\lvert}{\rvert}%
\DeclarePairedDelimiter\norm{\lVert}{\rVert}%

\makeatletter
\let\oldabs\abs
\def\abs{\@ifstar{\oldabs}{\oldabs*}}
\let\oldnorm\norm
\def\norm{\@ifstar{\oldnorm}{\oldnorm*}}
\makeatother

\DeclareMathSymbol{\sm}{\mathbin}{AMSa}{"39}